\documentclass[a4paper,11pt]{ijamas}

\newcommand{\pr}{\mbox{\sf P}}
\newcommand{\ex}{{\bf\sf E}}               %% expectation
\newcommand{\var}{\mbox{\sf Var}}

\newcommand{\call}{{\cal L}}

\newcommand{\calf}{{\cal F}}

\newcommand{\al}{\alpha}                %%
\newcommand{\g}{\lambda}               %%
\newtheorem{thm}{Theorem}
\newtheorem{lem}[thm]{Lemma}

\def\th{\theta}
\def\Th{\Theta}

\begin{document}

\baselineskip 18pt
\title{Approximating the Value Functions of Stochastic
Knapsack Problems: A Homogeneous Monge-Amp\'ere Equation and Its
Stochastic Counterparts}

\date{Yingdong Lu\\
IBM T.J. Watson Research Center\\ Yorktown Heights, NY 10598 }
\maketitle

\begin{abstract}
Stochastic knapsack problem originally was a versatile model for
controls in telecommunication networks. Recently, it draws
attentions of revenue management community by serving as a basic
model for allocating resources over time. We develop approximation
schemes for knapsack problems in this paper, a system of nonlinear
but solvable partial differential equations and stochastic partial
differential equation are shown to be the limit of the process that
following the optimal solution of the stochastic knapsack problem.

\medskip

\noindent\textbf{Keywords:} stochastic knapsack problem, fluid
limit, diffusion approximation, Monge-Amp\'ere eqnation.

\medskip

\noindent\textbf{2000 Mathematics Subject Classification:} 49J20,
65C30.

\end{abstract}

\vskip .5cm

\section{Introduction}
\label{sec:problem}

Stochastic knapsack problem(a.k.a. dynamic knapsack problem) is a
typical mathematical model for sequential resource allocation, see,
e.g. \cite{prastacos}. Similar model was proposed also for queueing
control problems, see, e.g. \cite{altman}, \cite{hajek} and
\cite{lippman}. In the 1990's, this type of models were intensively
studied in the analysis and design of telecommunication networks,
see, e.g. \cite{rossyao} and \cite{keith}. Lately, the stochastic
knapsack problem drew considerable attentions in the operations
research community due to their applications in revenue (yield)
management, see, \cite{lly2003} and references therein. The solution
of the problem can generally obtained through dynamic programming.
However, these stochastic knapsack problems are usually only
components of a larger optimal control problem. For example, in many
revenue management problems, stochastic knapsack serves as a
realistic model for pricing activities. In order to achieve overall
financial objectives, a pricing model have to be tied with
manufacturing and supply chain management. In these cases,
constantly checking the look-up table generated by dynamic
programming will not be feasible. Efficient approximations in more
explicit forms are hence needed. In this paper, we intend to study
the first and second order approximations in concise form.

The methodology we follow is basically the fluid and diffusion
approximation approaches that derived from functional law of large
numbers and functional central limit theorems. These methods have
been applied very successfully in the study of queues and queueing
networks, see, e.g. \cite{chenyaobook}. Our study can be viewed as a
further enrichment to this body of research. In queueing context,
the limit process are usually functional of simple differential
equations and (reflected) linear differential equations. The limit
we obtain is a solution of Monge-Amp\'ere equation, one of the most
important nonlinear partial differential equations.  Monge-Amp\'ere
equation was linked to crucial quantities in geometry, was
extensively studied by researchers in various branches of
mathematics. We also hope that this could lead to further
interactions between related fields.

The specific setup of the problem under consideration is the
following, at time $t=0$, we have $W$ units of resource available,
at any time $t=1, \cdots, T$, a request(demand) arrives in the form
of a bivariate random vector, $(P_t, Q_t)$, where the two components
represent the unit offer price and the required quantity; at each
time $t$, after the realization of the request, we decide whether to
accept or reject it to maximize the average revenue collected by
time $T$. Furthermore, we assume that,

\begin{itemize}
\item
Suppose $(P_t,Q_t)$, for each $t$ (and i.i.d.\ across $t$) follows
a joint (discrete) distribution:
\begin{eqnarray}
\label{jointdistn} \pr [P_t=p_i, Q_t=q_j]:=\th_{ij}, \qquad
i=0,1,...,\ell; \: j=1,...,k .
\end{eqnarray}
Here, $p_0=0$, and $\th_0:=\sum_j \th_{0j}$ represents the
probability that there is no arrival (request/demand) in a given
period;
\item
No partial fulfillment is allowed.
\end{itemize}

The decision of accept/rejection at each time $t$, given that the
available resource level is $d$, is determined by the following
stochastic dynamic programming,
\begin{eqnarray}
\label{dprecursion}
V(t,d)&=&V(t+1,d)[\th_0+\Th (d)]\nonumber\\
&+&\sum_{i\neq 0}\sum_{j:q_j\le d}\th_{ij}\cdot
\max\{p_iq_j+V(t+1,d-q_j), V(t+1,d)\},
\end{eqnarray}
where
\begin{eqnarray}
\label{Th} \Th(d):=\sum_{i\neq 0}\sum_{j:q_j > d}\th_{ij}.
\end{eqnarray}
Clearly, the first term on the right hand side of
(\ref{dprecursion}) corresponds to the case of either no arrival
or the request size exceeds the available inventory; whereas each
term under the summation compares the two actions: accept (i.e.,
supply) the request, or reject it. If we supply the request, then
we earn the revenue $p_jq_j$, and proceed to the next period with
$q_j$ units less in the available inventory.

At the last period, we have
\begin{eqnarray}
\label{boundary} V(T,d)= \sum_{i\neq 0}\sum_{j:q_j\le d}\th_{ij}
p_iq_j,
\end{eqnarray}
since clearly the best action is to supply any possible requests
using the remaining inventory.

The optimal decision is hence, to accept the demand $(P_t, Q_t)$,
if it satisfies,
\begin{eqnarray}
\label{opt_threshold}
 P_tQ_t+V(t+1,d-Q_t)\ge V(t+1,d);
\end{eqnarray}
and reject it otherwise.

This problem arises from many applications involving sequentially
allocating limit amount of resources to multiple classes of demands,
hence has been widely studied. The existing research can be divided
into two groups, one focuses on the structural properties of the
value function, usually through Markov decision process and
modularity argument, see, e.g. \cite{altman}, \cite{hajek}; the
other focus on a variety of heuristics policies and their
performance, much attentions was drawn some applications in the
field of revenue management and dynamic pricing, which is just
another mechanism to implement the accept/rejection decision, see,
e.g. \cite{lly2003}. In this paper, we examine the asymptotic
behaviors of the value function. Our approach is similar to the
``fluid limits'' and ``diffusion limits'' methods in queueing
theory, where time and space are properly scaled to induce
asymptotic behaviors of the systems. To be more specific, we scale
both time and space, which is the amount we allow to allocate, by a
constant $n$, then let $n$ go to infinity, and examine the
asymptotic behaviors of the value function under this situation.
Equivalently, we consider that following two processes,
\begin{eqnarray*}
{\bar V}^n(t,d)= \frac{1}{n} V(nt, nd),\qquad {\hat V}^n(t,d)=
\frac{1}{\sqrt{n}} V(nt, nd),
\end{eqnarray*}
In the following two sections, we will exam the asymptotic
behavior of ${\bar V}$ as $n\rightarrow \infty$, we show that it
converges to the solution of a boundary value problem of a
homogeneous Monge-Amp\'ere equation, give the general solution to
the equation, and solve it for some special case; then we in order
to characterize the rate of the convergence and the asymptotic
random behaviors, we identify the stochastic differential equation
that the limit of ${\hat V}^n(t,d)$ satisfies.

When the knapsack is characterized by more than one characters, we
have a multi-dimensional stochastic knapsack problem. It is
discovered recently to be an important model in the study of
inventory management. Unlike the one dimensional case, the dynamic
programming for the multi-dimensional stochastic knapsack problem
is exponential with respect to the size of the problem.
Approximations can certainly play even more important roles. Our
fluid and diffusion approximation are carried over to this
important model. A system of differential equation needs to be
solved to get the fluid limit.

The reset of the paper will be organized as the following, in Sec.
\ref{sec:convergence}, we will establish the fluid limit through
probabilistic arguments, and present the solution to the
Monge-Amp\'ere equation; in Sec. \ref{sec:diff}, diffusion
approximations are established for the case of unit demand; then,
the model is extended to the case of multi-dimensional stochastic
knapsack problem, for both fluid and diffusion approximation.

\section{Convergence Results}
\label{sec:convergence}

To facilitate the analysis, we will take a de tour. Instead of the
value function, we study the stochastic process whose mean will
achieve the value function. Let us denote $X_{t,d}(s), s=t, t+1,
\cdots, T$ to be the reward collected at time $s$ following the
optimal policy while starting at time $t$ with $d$ units available.
Apparently, $V(t,d)= \ex[X_{t,d}(T)]$. Define,
\begin{eqnarray}
\label{scaling} V^n(t,d) = \frac{1}{n}X_{nt,nd}(nT).
\end{eqnarray}
We expect that a law of large number type result exists, so that
$V^n(t,d)$ can converge to $V(t,d)$ for each $(t,d)$. Moreover, we
hope to extract the dynamic it satisfies, thus give us the
relationship we expect from $V(t,d)$.

To prove the convergence, we need,
\begin{lem}
{\rm Let $V^n(t, d)$ be the random variable defined above, then,
\begin{eqnarray}
\label{order_estimate} \var[X_{nt, nd} (nT) ]= O(n), \hbox{as} \ \
n\rightarrow \infty.
\end{eqnarray}}
\end{lem}
{\bf Proof} We prove this lemma by induction. First, when $n=1$,
observe that when $d$ approach infinity, the optimal policy will
be just accept any arrivals, then the variance will be constant.
Therefore, we have a maximum for the variance for all $d$. Denote
$M$ to be twice that number, we want to show that for each $n$,
$$\frac{\var[X_{nt, nd} (nT) ]}{n} \le M, \forall n, d.$$
Suppose it holds up to $n$. Then for the case of $n+1$,
apparently, we need,
\begin{eqnarray*}
\frac{\var[X_{(n+1)t, (n+1)d} ((n+1)T) ]}{n+1} \le M, \forall d,
\end{eqnarray*}
To see this, conditional upon the arrivals from $(n+1)t+1$ to
$(n+1)t+(T-t)$, whose $\sigma$-algebra can be denoted as
$\calf_{(n+1)t+(T-t)}$ we have,
\begin{eqnarray*}
\var[X_{(n+1)t, (n+1)d} ((n+1)T)]&=&\ex[\var[X_{(n+1)t, (n+1)d}
((n+1)T)|\calf_{(n+1)t+(T-t)}]]\\ &&+\var[\ex[X_{(n+1)t, (n+1)d}
((n+1)T)|\calf_{(n+1)t+(T-t)}]],
\end{eqnarray*}
in which,
\begin{eqnarray*}
\ex[\var[X_{(n+1)t, (n+1)d} ((n+1)T)|\calf_{(n+1)t+(T-t)}]]\le nM,
\end{eqnarray*}
from inductive assumption, and
\begin{eqnarray*}
\var[\ex[X_{(n+1)t, (n+1)d} ((n+1)T)|\calf_{(n+1)t+(T-t)}]]\le M.
\end{eqnarray*}
due to the definition of $M$. $\Box$

For the convergence proof, we adapt the classical proof of the
strong law of large number, see, e.g. \cite{billingsley}. Equation
(\ref{order_estimate}) and Markov inequality imply,
\begin{eqnarray*}
\sum_{n=1}^\infty \pr\left[ \left| \frac{X_{u_nt, u_nd}
(u_nT)-\ex[X_{u_nt, u_nd} (u_nT)]}{u_n}\right|>\epsilon \right]
<\infty,
\end{eqnarray*}
where $u_n=\lfloor \al^n \rfloor$ for some $\al>1$. From
Borel-Cantelli lemma, we know that $X_{u_nt, u_nd} (u_nT)/u_n$
converges almost surely, and let us denote the limit as $u(t,d)$ for
each $(t,d)$. Meanwhile, for any $u_n \le k \le u_{n+1}$, we have,
\begin{eqnarray*}
\frac{u_n}{u_{n+1}} \frac{X_{u_nt, u_nd} (u_nT)}{u_n} \le
\frac{X_{kt, kd} (kT)}{k} \le \frac{u_{n+1}}{u_n} \frac{X_{u_{n+1}t,
u_{n+1}d} (u_{n+1}T)}{u_{n+1}}.
\end{eqnarray*}
This leads to,
\begin{eqnarray*}
\frac{1}{\al} u(t,d) \le \lim \inf_k \frac{X_{kt, kd} (kT)}{k} \le
\lim \sup_k \frac{X_{kt, kd} (kT)}{k} \le \al u(t,d),
\end{eqnarray*}
since it holds for any $\al>1$, we can conclude that the
convergence with probability one.

Next, we relax the integral assumptions on the time and the
knapsack size. For any $(t, d)\in \mathbb{R}^2_+$, define, ${\hat
X}_{nt, nd} (nT)= X_{\lfloor nt\rfloor,\lfloor nt\rfloor} (nT)$.
Then it is easy to see that the convergence can be extended.
Furthermore, the converge is uniform on compact set, for related
details, see. e.g. \cite{chenyaobook}. Therefore we have,

\begin{thm}
\label{thm:slln} $V^n(t,d)\rightarrow u(t,d)$ almost surely as $n
\rightarrow \infty$. and $u(t,d)$ satisfies the Monge-Amp\'ere
equation with proper boundary value.
\begin{eqnarray}
\label{pde}
\left\{
 \begin{array}{c}\frac{\partial u(x,y)}{\partial x} +  g(
 \frac{\partial u(x,y)}{\partial y}) =0, \\
 u(X,y)=h(y), u(x,0)=0, (x,y) \in [0,X]\times[0, Y]
 \end{array}
 \right.
\end{eqnarray}
where $g(\cdot)$ is the loss function of the random variable $PQ$,
i.e. $g(x):=\ex[Q[P-x]^+]$.
\end{thm}
{\bf Proof} The only thing left to proof is that the limit will
satisfies the equation in (\ref{pde}). Since the convergence is
uniform on compact set, we can further conclude that, if $u$ is
continuously differentiable with respect to $t$ and $d$ on a compact
set, which imply uniform continuity for both $u$ and its
derivatives, then,\footnote{The following relations can be also
derived through actions of Schwartz functions to define derivative
in the distribution sense, hence, the solutions of the partial
differential differential equation are weak solutions. However, in
our case, the PDE has a classical solution, this direct approach is
more intuitive.}
\begin{eqnarray*}
X_{nd,nt}(nT) -  X_{nd,nt+1}(nT)=\frac{\frac{X_{nd,nt}(nT)}{n} -
\frac{X_{nd,nt+1}(nT)}{n}}{1/n}\rightarrow - \frac{\partial
u(t,d)}{\partial t}, {\rm a.s.}
\end{eqnarray*}
\begin{eqnarray*}
X_{nd,nt}(nT) -  X_{nd+1,nt}(nT)=\frac{\frac{X_{nd,nt}(nT)}{n} -
\frac{X_{nd+1,nt}(nT)}{n}}{1/n}\rightarrow - \frac{\partial
u(t,d)}{\partial d}. {\rm a.s.}
\end{eqnarray*}

This can be achieve from the following relationship.
\begin{eqnarray}
\label{difference} X_{nd,nt}(nT) -  X_{nd,nt+1}(nT) = [V(nd-Q,
nt+1)+PQ -V(nd, nt+1)]^+
\end{eqnarray}
This is the same as,
\begin{eqnarray}
\label{difference} \frac{\frac{X_{nd,nt}(nT)}{n} -
\frac{X_{nd,nt+1}(nT)}{n}}{\frac{1}{n}} = [PQ+\frac{\frac{V(nd-Q,
nt+1)}{n} -\frac{V(nd, nt+1)}{n}}{\frac{1}{n}}]^+.
\end{eqnarray}
Since, both $\frac{X_{nd,nt}(nT)}{n}$ and ${\bar V}^n(d,t)$, as its
mean, converge to the same function uniformly on  compact set, and
according to regularity theorems on the first order partial
differential equations, both this function and its derivatives are
uniformly continuous, see, e.g. \cite{john}, we can ensure the
convergence of differences to the derivatives. (\ref{difference})
gives us the relation,
$$\frac{\partial
u(x,y)}{\partial x} +  g(\frac{\partial u(x,y)}{\partial y}) =0.
$$ We can draw the same conclusion from the boundary value. $\Box$

The equation (\ref{pde}) is, of course, a nonlinear partial
differential equation, which in general is extremely difficult to
analyze. However, we can reduce this problem to a very well known
equation. Just take derivative with respect to $x$, we have,
\begin{eqnarray*}
u_{xx}+ g'(u_y)u_{xy}=0,
\end{eqnarray*}
similarly, take derivative with respect to $y$, we have
\begin{eqnarray*}
u_{xy}+ g'(u_y)u_{yy}=0,
\end{eqnarray*}
multiply the first one by $u_{yy}$ and the second one by $u_{xy}$,
then take difference, we have, $u_{xx}u_{yy}-u_{xy}^2=0$. This is a
homogeneous Monge-Amp\'ere equation. In general, Monge-Amp\'ere
equation is an important mathematical subject in analysis and
geometry. The general solution to the homogeneous Monge-Amp\'ere
solutions was obtained independently in \cite{zhdanov} and
\cite{ME}. We will follow the basic setup and structure of the
latter. In Fairlie and Leznov(1995) \cite{ME}, a general
construction to the second class is obtained. Given any two
arbitrary $C^1$ function $R(\cdot)$ and $f(\cdot)$, a solution can
be uniquely determined. Here, the general procedure of the
construction gives us the basic forms that the solution will take,
$g(\cdot)$ and $h(\cdot)$ then will help us to determine the
solution completely. We will also reveal the relations between the
two sets of functions. The general solution to the homogeneous
Monge-Amp\'ere equation bears the following form,
$$ \frac{\partial u}{\partial x}= R(\xi) -\xi \frac{d R}{d\xi},
\frac{\partial u}{\partial y}=\frac{dR}{d\xi},$$
\begin{eqnarray}
\label{linear} y-x\xi=f(\xi).
\end{eqnarray}
For any continuous function $R(\xi)$ and $f(\xi)$. For our problem,
we know that $\xi+ g(R') + g'(R') =0$, $u(X,y)=h(y)$, these
conditions can uniquely determine the two functions $R(\xi)$ and
$f(\xi)$, hence, uniquely determine $u(x,y)$. For example, when the
batch distribution takes the following special form, $g(x)=e^{-\g
x}$, $R'(\xi) = -\frac{\log(\frac{\xi}{\g-1})}{\g}$, we know that
\begin{eqnarray*}
u(x,y)= \int_{x_0}^x \left[R(\xi) -\xi \frac{d R}{d
\xi}\right]dx+\int_{y_0}^y \frac{d R}{d\xi} dy,
\end{eqnarray*}
in conjunction with the relation,
\begin{eqnarray*}
\frac{d \xi}{dx} =\frac{\xi}{x-f'(\xi)}, \frac{d
\xi}{dy}=\frac{1}{f'(\xi)-x}.
\end{eqnarray*}
we can obtain $u(x,y)$ from routine differential equation procedures
for first order partial differential equation, see e.g. \cite{john}.

%\section{Multidimensional Stochastic Knapsack Problem}
%\label{sec:mskp}

%Now, let us turn to the problem of multidimensional stochastic
%knapsack problem. First, let us introduce the problem briefly.
%Now, the capacity of the knapsack is $m$ dimensional vector,
%$(W_1, W_2, \cdots, W_m)$, at each time $t$, the arrival is a
%$m+1$ dimensional vector of $(p, q_1, q_2, \cdots, q_m)$, where
%$(q_1, q_2, \cdots, q_m)$ represents the capacity required, and
%$p$ represents rewards collected if the

\section{Second Order Approximations}
\label{sec:diff}

While fluid limits capture the dynamics of the system evolution,
asymptotic characterization of the randomness of the system is also
desired. Diffusion approximations, results of various central limit
theorems and strong approximation theorems usually serve this
purpose. In addition, diffusion approximations can also provide an
estimation of the order of the system converging to the fluid
limits. The process understudy differs from most of those appears in
the queueing context is that the main process is a random walk with
time and state dependency, furthermore, it depends upon the
increments of another random walk. Standard functional central limit
theorem can not be applied directly. How ever in the case of unit
demand, we can overcome this barrier by making use of results
derived in \cite{anisimov}. In \cite{anisimov}, a large class of
switching system are studied. Our system turns out to be one that
``switching'' at every time period. To apply the results in
\cite{anisimov}, it is require that the evolution function of the
mean and variance are smooth, in fact Lipschitz should be enough.
Those functions for our case are functions of the solutions to the
Monge-Amp\'ere equation, whose smoothness has been well established.
Then, results in \cite{anisimov} will give us the stochastic
differential equation(SDE) the limiting process will satisfy.
However, the existence and uniqueness of the solution to the SDE
were not established in \cite{anisimov}. Meanwhile, in
\cite{krylov}, the existence and uniqueness of the solution to this
type of SDE with more general assumptions on the coefficients have
been obtained. Hence, the SDE uniquely determine a stochastic
process.

To capture the randomness of $X_{n,d}(T)$, we define $y_{d,T}(s)$
to be the number of units that are supplied to the demands up to
time $s$. From the set up, we know that $y_{d,T}(s)$ is a Markov
chain. Moreover, we have,
\begin{eqnarray*}
y_{d,T}(s+1)=y_{d,T}(s)+\chi(\{V(s+1, d-y_{d,T}(s)-1)+p-V(s+1,
d-y_{d,T}(s))>0\}),
\end{eqnarray*}
\begin{eqnarray*}
Z_{d,T}(s)-Z_{d,T}(s)= p_i[y_{d,T}(s)-y_{d,T}(s)].
\end{eqnarray*}

In \cite{anisimov}, random walks of this type have been studied.
More specifically, Let $\{\xi_n(\al),\al \in\mathbb{R}_+\}$ be
independent random variables with $al$ as parameters, and $S_{n+1} =
S_n+ \xi(S_n)$. The following theorem is a summary of results in
\cite{anisimov},
\begin{thm}{\rm (\cite{anisimov})}
If there exists a function $s(t)$ such that,
\begin{eqnarray*}
|\frac{1}{n}S_n(nt)-s(t)| \rightarrow 0, a.s.
\end{eqnarray*}
and $D_n\rightarrow D(\al)$, then, $D_n= \var[\xi_n(\al)]/n$ then,
$\gamma_n (t) := n^{-1/2} [S_n(nt)-ns(t)]$ weekly converge in $D$
to the following stochastic differential equation,
\begin{eqnarray*}
d \gamma(t) = D(s(t))^{1/2} dw_t,
\end{eqnarray*}
provided that its solution exists and is unique.
\end{thm}
Following the similar functional strong law of large numbers in
the previous section, we can see that $s(t)$ satisfies,
\begin{eqnarray*}
\frac{ds(t)}{t} =F(u_d(t, d-s(t))).
\end{eqnarray*}
By the Lipschitz continuity of the solution to the Monge-Amp\'ere
equation, and existence and stability theorems of ordinary
differential equation, we can see that $s(t)$ can be uniquely
determined and is Lipschtz continuous. Therefore, we can conclude
that $n^{-1/2} [Y_n(nt)-ns(t)]$ converge to a stochastic process
that satisfies the following stochastic differential equation,
\begin{eqnarray}
\label{sde1} dY_t=   \sqrt{g(u_d(t, d-s(t)))[1-g(u_d(t,
d-s(t)))]}dW_t,
\end{eqnarray}
where $dW_t$ denotes a wiener integral with respect to a standard
Brownian motion, and the revenue collected can be represented by
\begin{eqnarray}
\label{sde2} dZ_s= PdY_s.
\end{eqnarray}

\noindent {\bf Remark: } For the problems with general demand batch
size, we conjecture that the process will weakly converge to a
process that satisfies the following stochastic differential
equation,
\begin{eqnarray*}
dY_t= \sqrt{F_{p,d}(u_d(t, d-Y_t))[1-F_{p,d}(u_d(t, d-Y_t))]}dW_t.
\end{eqnarray*}
Where $F_{p,d}$ is the distribution function for both price and the
quantity, and is not necessary continuous. Hence, the limiting
process is not a diffusion process any more. However, results in
\cite{krylov} guarantee the existence and uniqueness of the limiting
process. This convergence will be the subject of future research.

\section{Multi-dimensional Stochastic Knapsack Problems}
\label{sec:multi}

In this section, we intend to deal with a multi-dimensional
version of the stochastic knapsack problem. To be more specific,
at time $t=0$, we have a vector of $m$ different types of products
available, $(W^1,W^2,\cdots, W^m)$, at any time $t=1,2,\cdots, T$,
the demand arrives now is characterized by a $m+1$ tube,
$(P_t,Q^1_t, Q^2_t,\cdots, Q^m_t)$, where $P_t$ describes the
rewards and a $m$-dimensional characterization vector $(Q^1_t,
Q^2_t,\cdots, Q^m_t)$ describes the combination of the the $m$
types products it requires. They are i.i.d random variables, with
discrete distribution,
\begin{eqnarray*}
\pr[P_t =p_i, Q^1_t=q_{j_1}^1, Q^2_t=q_{j_2}^2,\cdots,
Q^m_t=q_{j_m}^m] = \th_{i,j_1,j_2, \cdots, j_m}.
\end{eqnarray*}
Again, at each time $t$, we need to determine whether to accept or
to reject the demand arrival to achieve the maximum average
revenue at time $T$.  Let us denote $V(t, d^1, d^2,\cdots, d^m)$
be the value function for the dynamic programming. Similarly to
the one dimensional case, it satisfies the following recursive
equation,
\begin{eqnarray}
\label{dp_multi}&& V(t, d^1, d^2,\cdots,d^m)\nonumber\\ &=& V(t+1,
d^1, d^2,\cdots,d^m)[\th_0+\Th(d^1, d^2,\cdots,d^m)] \nonumber\\
&& + \sum \th_{i,j_1,j_2, \cdots, j_m} \max\left\{ p_i +V(t+1,
d^1-q_1, d^2-q_2,\cdots,d^m-q_m), \right. \nonumber \\ && \left.
V(t+1, d^1, d^2,\cdots,d^m)\right\}.
\end{eqnarray}
Here,
$$\Th(d^1, d^2,\cdots,d^m)=\pr[\cup_{k=1}^m\{ Q^k_t>d^k\}]. $$
It is obvious that the dynamic programming is not so easy-to-solve
as those in one-dimensional case. In fact, the computational
efforts needed grow exponentially with respect to $W=\max_k W^k$.
Therefore, approximation will be much more desirable.

Apply the sam scaling and argument, we can see that the similar
convergence results hold,
\begin{thm}
\label{thm:slln_multi} $V^n(t,d^1, d^2,\cdots,d^m)\rightarrow
u(t,d^1, d^2,\cdots,d^m)$ as $n \rightarrow \infty$. and $u(t,d^1,
d^2,\cdots,d^m)$ satisfies the following boundary value problem,
\begin{eqnarray*}
\left\{
 \begin{array}{c}\frac{\partial u(x_0, x_1, x_2, \cdots,x_{n-1})}{\partial x_0} +  G(
 \frac{\partial u}{\partial x_1}, \frac{\partial u}{\partial x_2},\cdots,\frac{\partial u}{\partial x_{m-1}}) =0, \\
 u(x_0,x_1, x_2, \cdots, x_{m-1})=h(x_1, x_2, \cdots, x_{n-1}), u(x_0,0)=0,
 \end{array}
 \right.
\end{eqnarray*}
where
$$G(z_1,z_2,\cdots,z_{m-1}):=\ex[P-(z_1+z_2+\cdots+z_{m-1})]^+.$$
Hence, $u$ satisfies the homogeneous Monge-Amp\'ere equation
$Det(D^2u)=0$.
\end{thm}
{\bf Proof: } The only thing left is to show that $u$ satisfy the
Monge-Amp\'ere equation. To see this, take derivative with respect
to $x_0, x_1, \cdots, x_m $ respectively on equation,
$$ \frac{\partial u(x_0, x_1, x_2, \cdots,x_{n-1})}{\partial x_0} +  G(
 \frac{\partial u}{\partial x_1}, \frac{\partial u}{\partial x_2},\cdots,\frac{\partial u}{\partial x_{m-1}})
 =0, $$
we have, for $k=0,1,\cdots, m$

\begin{eqnarray}
\label{linear}\frac{\partial u^2}{\partial x_0x_k} +
\sum_{i=1}^{m-1}G_i(
 \frac{\partial u}{\partial x_1}, \frac{\partial u}{\partial x_2},\cdots,
 \frac{\partial u}{\partial x_{m-1}})\frac{\partial u^2}{\partial x_ix_k}=0,
\end{eqnarray}
regard this as a linear system for variables $(1, G_1, G_2, \cdots,
G_m)$, then in order for (\ref{linear}) to hold, the coefficient
matrix, $(D^2u)$, must be singular. $\Box$

Again from \cite{ME}, we have the general solution for $n$
dimensional Monge-Amp\'ere equation. Understandably, it depends
upon the selection of two $m-1$ variate functions.  Given
arbitrary function $R(\xi_1, \xi_2, \cdots, \xi_{n-1})$, $L(\xi_1,
\xi_2, \cdots, \xi_{n-1})$,
$$u_0=R-\sum\xi_jR_j, u_j =R_j,$$
$$ x_j-\xi_j x_0= Q^j(\xi_1, \xi_2, \cdots, \xi_{n-1}),$$
$$Q=(D^2R)^{-1}DL.$$
Where $DL$ and $D^2R$ denote the gradient vector and Hessian
matrix of L and R respectively. Function $G(\cdot)$ and the
boundary condition $h(\cdot)$ will then determine $R$ and $L$.

Define $y^k_{d,T}(s), k=1,2,\cdots, m$ to be the number of units
of characterizations that are supplied to the demands up to time
$s$. There exists a function $(s^1(t), s^2(t),\cdots, s^m(t))$
which is the solution to the following differential equation
system,
\begin{eqnarray*}
\frac{ds^k(t)}{t} =G(u_i(t, d^1-s^1(t),d^2-s^2(t) ,\cdots,
d^m-s^m(t)), \qquad k=1,2,\cdots, m.
\end{eqnarray*}
The regularity of the solution will then help us conclude that,
$$\left(\frac{Y^k_n(nt)-ns^k(t)}{n^{1/2}}\right),$$
converges to a stochastic process determined by the following
stochastic differential equation,
\begin{eqnarray}
\label{sde2} dY^k_t=   \sqrt{{\tilde G_k}(1-{\tilde G_k})}dW^k_t,
\end{eqnarray}
where, ${\tilde G_k}=G(u_k(t, d^1-s^1(t),d^2-s^2(t) ,\cdots,
d^m-s^m(t))$.

\bibliographystyle{dcu}
\bibliography{IJAMAS_example}

\end{document}